%%PlainTeX

\tolerance=8000
\def\detail#1{}
%---------
\def\redef#1#2{\expandafter\ifx\csname #1\endcsname\relax
    \expandafter\edef\csname #1\endcsname{#2}
    \else \message{redefinition of \string#1} \hiba \fi }

\openin 0=cite.inc
\ifeof 0
\closein 0
\message{One more pass needed for references }
\else
\closein 0
\input cite.inc
\fi

\immediate\openout 0=cite.inc
\immediate\openout 1=pages.txt

% table of contents
\immediate\openout 2=conten.inc
\def\writeitem#1#2{\write2
  {\string\line{#1\string\nagy{}#2
  \string\rm\string\leaderfill\string\quad\folio}}}

\def\nagy{\font\caps=cmcsc10 \caps}

\magnification=\magstep1

\hoffset=0.5 true in
\hsize=6.0 true in
\voffset=.3 true in
\vsize=8.28 true in

\def\boxit#1{\vbox{\hrule\hbox{\vrule#1\vrule}\hrule}}
\def\items#1{\itemitem{{\rm #1}}}

\newcount\sorszam \sorszam=0
\def\sorszaminc {\advance\sorszam by1 \the\section.\the\sorszam. }

\def\uj {\bigskip \rm}

\def\oldalszam{
\headline={
 \nagy
 \ifnum \folio > 1
 \ifodd \folio {\hfil n\'andor sieben \hfil} \else
 {\hfil crossed products by semigroup actions \hfil}
 \fi
 \folio
% \write1{\the\inputlineno}
 \else \hfil
 \fi
}
}

\newcount\section \section=0
\def\newsection#1{
%     \vfil\eject
     \goodbreak
     \advance\section by1 \sorszam=0
     \medskip\bigskip\centerline{\nagy \the\section. #1}
% for table of contents
     \writeitem{\the\section.  }{#1}
     \nobreak
     }

\def\definition {\bigskip\goodbreak\noindent\bf Definition \sorszaminc \rm \ }
\def\theorem {\bigskip\goodbreak\noindent \bf Theorem \sorszaminc \sl \ }
\def\prop {\bigskip\goodbreak\noindent \bf Proposition \sorszaminc \sl \ }
\def\lemma {\bigskip\goodbreak\noindent \bf Lemma \sorszaminc \sl \ }
\def\cor  {\bigskip\goodbreak\noindent \bf Corollary \sorszaminc \sl \ }
\def\remark {\bigskip\goodbreak\noindent \bf Remark \sorszaminc \rm \ }
\def\proof { \medskip {\noindent \it Proof. } \rm }
\def\example {\bigskip\goodbreak\noindent \bf Example \sorszaminc \rm \ }

\def\nullbox{\setbox0=\null \ht0=5pt \wd0=5pt \dp0=0pt \box0}
\def\eop {\nobreak\hfill \boxit{\nullbox}\goodbreak}
\def\Eop {\nobreak\ \hfill \boxit{\nullbox}}

\def\hivat#1{  %  -#1-
\immediate
\write0{\string\redef{\string#1}{{\string\rm\the\section.\the\sorszam}}}}

\def\hivatt#1{ %\edef#1{{\rm\the\section.\the\sorszam}}
\immediate
\write0{\string\redef{\string#1}{{\string\rm\the\sorszam}}}}
\def\cite#1{\csname #1\endcsname}

\def\a {\ifmmode \alpha \else $\alpha$ \fi}
\def\b {\ifmmode \beta \else $\beta$ \fi}
\def\d {\ifmmode \delta \else $\delta$ \fi}
\def\PAut{\hbox{\rm PAut}}
\def\CovRep{\hbox{\rm CovRep}}
\def\slim{\mathop{\hbox{\rm s-lim}}}

\vsize=8.65 true in

\footline={\tenrm \hfil \folio \hfill}
\footline={}
\oldalszam

%%%%%%%%%%%%%%%%%%%%%%%%%%%%%%%%%%%%%%%%%%%%%%%%%%%%%%%%%%%%%
% \double
\bigskip
\centerline{$C^*$-CROSSED PRODUCTS BY PARTIAL ACTIONS AND}

\centerline{ACTIONS OF INVERSE SEMIGROUPS}
\uj
\centerline{N\'andor Sieben}

\footnote{}{{\it Date. }\the\month/\the\day/\the\year


This material is based upon work supported
by the National Science Foundation under Grant No.
DMS9401253.}

\centerline{}
\bigskip

{
\narrower\narrower
\centerline{\nagy Abstract}
\uj
%\sevenrm
%\baselineskip=.7\baselineskip
The recently developed theory of partial actions of discrete groups on
$C^*$-algebras is extended.
A related concept of actions of inverse semigroups
on $C^*$-algebras is defined, including covariant representations and
crossed products. The main result is that every partial crossed product
is a crossed product by a semigroup action.

}

\newsection{Introduction}
\uj
Recently the notion of a partial crossed product of a $C^*$-algebra by a
discrete group was defined by McClanahan [\cite{mcc}] as a generalization
of Exel's definition in [\cite{ex}].
The more well-established notion of the
crossed product of a $C^*$-algebra by an action of a group
uses a homomorphism into the automorphism group of
the $C^*$-algebra. The idea of a partial action is to replace the
automorphism group by the inverse semigroup of partial automorphisms.
A partial automorphism is an isomorphism between two closed ideals of a
$C^*$-algebra. Of course
we cannot talk about a homomorphism from a group into an inverse semigroup;
a partial action is an appropriate generalization.
In Section~2 we give a detailed discussion of partial actions.

After replacing the automorphism group by the semigroup of partial
automorphisms the next natural step is to replace our group by an inverse
semigroup. This makes it possible to use a more natural semigroup homomorphism
instead of a partial action.
We develop the elementary theory of an action of an inverse semigroup
in Section~4, and
we define the crossed product
by an inverse semigroup action in Section~5.

It turns out that there is a close connection between partial crossed products
and crossed products by inverse semigroup actions. In Section~6 we explore
this connection, showing that every partial crossed product is isomorphic to
a crossed product by an inverse semigroup action.

The majority of the research for this paper was carried
out while the author was a student at Arizona State
University. The results formed the author's masters
thesis written under the supervision of John Quigg.

I would like to take the opportunity to thank Professor
Quigg for his help and guidance during the writing of
this thesis.
% \single

\newsection{Partial actions}

\uj
In this section we discuss the notion of {\it partial actions} defined
by McClanahan [\cite{mcc}] which is a generalization of Exel's definition in
[\cite{ex}]. The major new result is Theorem \cite{hive}.

\definition
Let $A$ be a $C^*$-algebra. A {\it partial automorphism } of $A$ is a triple
$(\a,I,J)$ where $I$ and $J$ are closed ideals in $A$ and $\a:I
\rightarrow J$ is a *-isomorphism. We are going to write $\a$ instead of
$(\a,I,J)$ if the domain and range of $\a$ are not important.

\uj
If $(\a,I,J)$ and $(\b,K,L)$ are partial automorphisms of $A$
then the product $\a\b$ is defined as the composition of \a and \b
with the largest possible domain, that is,
$\a\b:\b^{-1}(I)\rightarrow A$, $\a\b(a)=\a(\b(a))$.
It is clear that $\b^{-1}(I)$ is a closed ideal of $K$. Since a closed ideal
of a closed ideal of $A$ is also a closed ideal of $A$, the product
$(\a\b,\b^{-1}(I),\a\b(\b^{-1}(I)))$ is a partial automorphism too.

A semigroup $S$ is an {\it inverse semigroup} if for every $s\in S$ there
exists a unique element
$s^*$ of $S$ so that $ss^*s=s$ and $s^*ss^*=s^*$.
The map $s\mapsto s^*$ is an involution.
An element $f\in S$ satisfying $f^2=f$ is called an {\it idempotent} of $S$,
and in this case $f=f^*$ also.
The set of idempotents of an inverse semigroup is a semilattice
with $e\wedge f=ef$.
Our general reference on semigroups is [\cite{how}].
It is easy to see that the set $\PAut(A)$ of partial automorphisms of $A$
is a unital inverse semigroup with identity $(\iota,A,A)$, where $\iota$ is the
identity map on $A$, and $(\a,I,J)^*=(\a^{-1},J,I)$.

\definition
Let $A$ be a $C^*$-algebra and $G$ be a discrete group with identity $e$.
A {\it partial action} of
$G$ on $A$ is a collection $\{ (\a_s,D_{s^{-1}},D_s) : s\in G \}$ of
partial automorphisms (denoted by \a or by $(A,G,\a)$ ) such that
\items{(i)} $D_e=A$
\items{(ii)} $\alpha_{st}$ extends $\alpha_s\alpha_t$, that is,
$\a_{st}|\a_t^{-1}(D_{s^{-1}})=\a_s\a_t$
for all $s,t\in G$.
\hivat{hiv77}

\prop
If \a is a partial action of $G$ on $A$ then
\items{(i)} $\a_e$ is the identity map $\iota$ on $A$
\items{(ii)} $\a_{s^{-1}}=\a_s^{-1}$ for all $s\in G $.
\hivat{hiva}

\proof
The statements follow from the following two identities
$$ \iota=\a_e\a_e^{-1}=\a_{ee}\a_e^{-1}=\a_e \a_e \a_e^{-1}=\a_e \  ,$$
$$\a_s\a_{s^{-1}}=\a_{ss^{-1}}|D_s=\a_e|D_s=\iota|D_s \  . $$
\Eop

\detail{
\uj
If $I$ and $J$ are ideals of $A$ then let $IJ$ denote the closure of the set
$\{ab:a\in I, b\in J\}$. Clearly $IJ \subset I \cap J$. On the other hand
if $c\in I\cap J$ and $e_\lambda$ is an approximate identity for $J$ then
$c=\lim_\lambda c e_\lambda$. This means that $c$ is in $IJ$ and so
$I\cap J=IJ$.
}

\lemma
If \a is a partial action of $G$ on $A$ then
$\a_t(D_{t^{-1}}D_s)=D_tD_{ts}$ for all $s,t \in G$.
\hivat{hivs}

\proof
By Proposition \cite{hiva} (ii),
$
\a_t(D_{t^{-1}}D_s)=\a_{t^{-1}}^{-1}(D_{t^{-1}}D_s)=\a_{t^{-1}}^{-1}(D_s)
$
which is the domain of $\a_{s^{-1}}\a_{t^{-1}}$ and hence is contained in the
domain $D_{(s^{-1}t^{-1})^{-1}}=D_{ts}$ of $\a_{s^{-1}t^{-1}}$. Since the
range of $\a_t$ is $D_t$,  we have
$$
\a_t(D_{t^{-1}}D_s)\subset D_tD_{ts} \hbox{\qquad for all \quad} s,t\in G.
$$
Since $\a_t$ is an isomorphism, the containment above implies
$$
D_{t^{-1}}D_s\subset \a_t^{-1}(D_tD_{ts})=\a_{t^{-1}}(D_tD_{ts})
  \hbox{\qquad for all \quad} s,t\in G .
$$
Replacing $t$ by $t^{-1}$ and $s$ by $ts$ gives
$$D_tD_{ts}\subset \a_t(D_{t^{-1}}D_{t^{-1}ts})=\a_t(D_{t^{-1}}D_s)
\hbox{\qquad for all \quad} s,t\in G $$
as desired.
\eop

\lemma
If \a is a partial action of $G$ on $A$ then
$\a_t(D_{t^{-1}}D_{s_1}\cdots D_{s_n})=D_tD_{ts_1}\cdots D_{ts_n}$ for all
$t,s_1,\ldots ,s_n \in G$.
\hivat{hivf}

\proof
The statement follows from the following calculation using Lemma \cite{hivs}:
$$\eqalign{
\alpha_t (D_{t^{-1}} D_{s_1}\cdots D_{s_n})
& = \a_t(D_{t^{-1}}D_{s_1}\cdots D_{t^{-1}}D_{s_n}) \cr
& = \a_t(D_{t^{-1}}D_{s_1})\cap\cdots\cap\a_t(D_{t^{-1}}D_{s_n}) \cr
& = D_tD_{ts_1}\cap\cdots\cap D_tD_{ts_n} \cr
& = D_tD_{ts_1}\cdots D_{ts_n} .\cr
} $$
\eop

\theorem
If \a is a partial action of $G$ on $A$ then the partial automorphism
$\a_{s_1} \cdots \a_{s_n}$ has domain
$D_{s_n^{-1}}D_{s_n^{-1}s_{n-1}^{-1}} \cdots D_{s_n^{-1} \cdots s_1^{-1} }$
and range
$D_{s_1}D_{s_1s_2} \cdots D_{s_1 \cdots s_n }$
for all $s_1,\ldots,s_n \in G$.
\hivat{hivc}
\hivat{hive}

\proof
The theorem is proven inductively. The statement is clear for $n=1$.
The induction step follows from the following calculation using Lemma
\cite{hivf}:
$$ \eqalign{
\hbox{ domain } \a_{s_1} \cdots \a_{s_n}
& = \a_{s_n}^{-1}(\hbox{ domain } \a_{s_1} \cdots \a_{s_{n-1}}) \cr
& = \a_{s_n^{-1}}
  (D_{s_n}D_{s_{n-1}^{-1}} \cdots D_{s_{n-1}^{-1} \cdots s_1^{-1} }) \cr
& = D_{s_n^{-1}}D_{s_n^{-1}s_{n-1}^{-1}} \cdots
   D_{s_n^{-1} \cdots s_1^{-1} } .
}
$$
The other part now follows since the range of
$\a_{s_1} \cdots \a_{s_n}$
is the domain of
$\a_{s_n^{-1}}\cdots \a_{s_1^{-1}}$.
\eop

\cor
The conditions in the definition of a partial action can be reformulated
as follows: \nobreak
\items{(i) } $D_e=A$
\items{(ii)$^\prime$} $\a_{st}|D_{t^{-1}}D_{t^{-1}s^{-1}}=\a_s\a_t$.

\remark
If \a is a partial action of $G$ on $A$ then \a generates a unital inverse
subsemigroup
$S=\{\a_{s_1} \cdots \a_{s_n} : s_1,\ldots, s_n \in G, n\in {\bf Z} \}$
of the semigroup of partial automorphisms of $A$.
Theorem \cite{hive} tells us the domains and ranges.
\hivat{hivm}

\newsection{Covariant representation}

\uj
We continue our discussion of partial actions. Here the major new results
are Theorem \cite{hiv1} and Corollary \cite{hivq}.

\definition
Let \a be a partial action of $G$ on $A$. A {\it covariant representation }
of \a is a triple $(\pi,u,H)$ where $\pi: A\rightarrow  B(H)$ is
a nondegenerate representation of $A$ on the Hilbert space $H$ and
$g\mapsto u_g:G \rightarrow B(H)$, where $u_g$ is a partial isometry
on $H$ with initial space $\pi(D_{g^{-1}})H$ and final space $\pi(D_g)H$,
such that
\items{(i)} $u_g\pi(a)u_{g^{-1}}=\pi(\a_g(a))\qquad$ for all
  $\quad a\in D_{g^{-1}}$
\items{(ii)} $u_{st}h=u_s u_t h \qquad$  for all
  $\quad h\in \pi(D_{t^{-1}}D_{t^{-1}s^{-1}})H$.
\hivat{hivd}

\uj
Notice that by the Cohen-Hewitt factorization theorem [\cite{coh}],
$\pi(D_g)H$
is a closed subspace of $H$ and so the notations for the initial and
final spaces make sense.

\prop If $(\pi,u,H)$ is a covariant representation then $u_e=1_H$ (the
identity map on $H$) and
$u_{s^{-1}}=u_s^*$ for all $s\in G$.
\hivat{hivn}

\proof
Since $\pi$ is a nondegenerate
representation, $H$ is the closed span of $\pi(A)H=\pi(D_e)H$.
But since $\pi(D_e)H$ is a closed subspace of $H$ we have $H=\pi(D_e)H$.
Hence $u_e$ has
initial and final space $H$ which means $u_e$ is a unitary on $H$.
By Definition \cite{hivd}\ (ii) we have $u_{e}=u_e u_e$ which implies
$1_H=u_e u_e^{-1}=u_e u_e u_e^{-1}=u_e$.

For the second statement we have to show that
$\langle u_s h,k\rangle=\langle h,u_{s^{-1}}k\rangle$ for all $h,k\in H$.
We can write $h=h_1+h_2$ where $h_1\in \pi(D_{s^{-1}})H$,
$h_2\in (\pi(D_{s^{-1}})H)^\perp$
and $k=k_1+k_2$ where $k_1\in \pi(D_s)H$,
$k_2\in (\pi(D_s)H)^\perp$. So it suffices to show that
$$\langle u_s h_1,k_1+k_2\rangle=\langle h_1+h_2,u_{s^{-1}}k_1 \rangle .$$
Since $\langle u_s h_1,k_2 \rangle=0=\langle h_2,u_{s{-1}} k_1 \rangle$, it
suffices to show that
$$\langle u_s h_1,k_1\rangle = \langle h_1,u_{s^{-1}}k_1 \rangle
\qquad \hbox{for all}\quad h_1\in \pi(D_{s{-1}})H, \quad k_1\in \pi(D_s)H .$$
Since $h_1=u_{s^{-1}}l$ for some $l\in \pi(D_s)H$ it remains to show that
$$\langle u_s u_{s^{-1}}l,k_1\rangle =
            \langle u_{s^{-1}}l,u_{s^{-1}}k_1 \rangle
  \qquad \hbox{for all}\quad k_1, l\in \pi(D_s)H .$$

By Definition \cite{hivd}\ (ii) we have
$u_em=u_s u_{s^{-1}}m$ for all $m \in \pi(D_s)H$. Hence $u_s u_{s^{-1}}l=l$
and the statement follows from the fact that $u_{s^{-1}}$ is a partial isometry
with initial space $\pi(D_s)H$.
\eop

\uj
Let $\pi_u :A\rightarrow B(H_u)$ be the universal representation of a
$C^*$-algebra $A$ on the
universal Hilbert space $H_u$. If $I$ is an ideal of $A$ then the double dual
$I^{**}$ of $I$, identified with the strong operator closure of $\pi_u(I)$, is
an ideal of the enveloping von
Neumann algebra $A^{**}$ of $A$, which is identified with the strong operator
closure of $\pi_u(A)$.
\detail{To see this let $i=\lim_\lambda i_\lambda$ be a strong operator limit
where
$i_\lambda\in\pi_u(I)$ and $a=\lim_\mu a_\mu$ be a strong operator limit
where $a_\mu\in\pi_u(A)$. Then $ia=(\lim_\lambda i_\lambda)(\lim_\mu a_\mu)=
\lim_\mu(\lim_\lambda i_\lambda a_\mu)$. Since $I$ is an ideal,
$i_\lambda a_\mu$ is an element of $\pi_u(I)$ and so
$\lim_\lambda i_\lambda a_\mu$ is an element of $I^{**}$ and hence $ia$ is
also an element of $I^{**}$. Similarly $ai$ is in $I^{**}$.
}
As such, $I^{**}$ has the form $pA^{**}$ for some central projection $p$
in $A^{**}$.

\definition
Let \a be a partial action of $G$ on $A$.
For $s\in G$, $p_s$
denotes the central projection of $A^{**}$ which is the identity of $D_s^{**}$.

\uj
Let $(\pi,u,H)$ be a covariant representation of $(A,G,\a)$.
Since $\pi$ is a nondegenerate
representation of $A$, $\pi$ can be extended to a normal morphism of
$A^{**}$ onto $\pi(A)^{\prime\prime}$. We will denote
this extension also by $\pi$.

\lemma
Let $\pi$ be a representation of $A$ on  $H$, $I$ be a closed ideal of
$A$ and $p$ be the central projection of $A^{**}$ which is the identity of
$I^{**}$. Then $\pi(I)H=\pi(p)H$.
\hivat{hivh}

\proof
If $a\in I$ and $h\in H$
then $\pi(a)h=\pi(pa)h=\pi(p)\pi(a)h\in\pi(p)H$
which implies $\pi(I)H\subset \pi(p)H$.

On the other hand $\pi(p)$ is in the strong operator closure of $\pi(I)$
and hence $\pi(p)$ is the strong operator limit of a net
$\{\pi(a_\lambda)\}$ in $\pi(I)$. Hence $\|\pi(p)h-\pi(a_\lambda)h\|
\rightarrow 0$ for all $h\in H$. Since $\pi(I)H$ is closed, this means
$\pi(p)h\in\pi(I)H$ for all $h\in H$.
\eop

\cor
Let $\pi$ be a representation of $A$ on  $H$. Then
$$\pi (D_{s_1}\cdots D_{s_n})H=\pi (p_{s_1}\cdots p_{s_n})H$$
for all $s_1,\ldots ,s_n\in G$.
\hivat{hivi}

\proof
It is clear that $p_{s_1}\cdots p_{s_n}$ is the identity of
$(D_{s_1}\cdots D_{s_n})^{**}$ in $A^{**}$, so the statement follows from
Lemma \cite{hivh}.
\eop

\cor
\hivat{helo}
If $(\pi,u,H)$ is a covariant representation then $u_su_s^*=\pi(p_s)$
and $u_s^*u_s=\pi(p_{s^{-1}})$ for all $s\in G$.
\eop

\prop
\hivat{hpppp}
If $(\pi,u,H)$ be a covariant representation then for $g_1,\ldots,g_n\in G$
we have
$$\eqalign{
u_{g_1}\cdots u_{g_n}u_{g_n}^*\cdots u_{g_1}^*&=
 \pi(p_{g_1}\cdots p_{g_1\cdots g_n}) \cr
u_{g_n}^*\cdots u_{g_1}^*u_{g_1}\cdots u_{g_n}&=
 \pi(p_{g_n^{-1}}\cdots p_{g_n^{-1}\cdots g_1^{-1}}) \cr
}$$

\proof
The second equality follows from the first one taking conjugates. For $n=1$
the first equality is true by Corollary \cite{helo}. Applying induction we
get
$$\eqalign{
u_{g_1}\cdots u_{g_n}u_{g_n}^*\cdots u_{g_1}^*&=
u_{g_1}u_{g_1}^*u_{g_1}\circ\pi(p_{g_2}\cdots p_{g_2\cdots g_n})\circ u_{g_1}^*
\cr
&=u_{g_1}\circ\pi(p_{g_1^{-1}}p_{g_2}\cdots p_{g_2\cdots g_n})\circ u_{g_1}^*
\cr
&=\pi(\a_{g_1}(p_{g_1^{-1}}p_{g_2}\cdots p_{g_2\cdots g_n}))\cr
&=\pi(p_{g_1}\cdots p_{g_1\cdots g_n}),\cr
}$$
by Lemma \cite{hivf}.
Notice that we extended $\a_{g_1}$ to get an isomorphism
$\a_{g_1}:D_{s^{-1}}^{**}\to D_{s}^{**}$ between the double duals.
\eop

\theorem
Let $(\pi,u,H)$ be a covariant representation. Then $u_{s_1}\cdots u_{s_n}$ is
a partial isometry with initial and final spaces
$$\pi (D_{s_n^{-1}}D_{s_n^{-1}s_{n-1}^{-1}} \cdots
  D_{s_n^{-1} \cdots s_1^{-1} })H \quad \hbox{and} \quad
\pi(D_{s_1}D_{s_1s_2}\cdots D_{s_1\cdots s_n})H, $$
for all $s_1,\ldots,s_n\in G$.
\hivat{hiv1}

\proof
It is clear from Proposition \cite{hpppp} that $u_{s_1}\cdots u_{s_n}$ is a
partial isometry. The statement about the initial and final spaces follows
from Corollary \cite{hivi}.
\eop

\cor
If $(\pi,u,H)$ is a covariant representation then
$$
u_{s_1\cdots s_n}h=u_{s_1}\cdots u_{s_n}h \quad \hbox{for all}\quad
h\in \pi (D_{s_n^{-1}}D_{s_n^{-1}s_{n-1}^{-1}} \cdots
  D_{s_n^{-1} \cdots s_1^{-1} })H
$$
and
$$\pi(a)u_{s_1\cdots s_n}=\pi(a)u_{s_1}\cdots u_{s_n}
\quad \hbox{for all} \quad a\in D_{s_1}D_{s_1s_2}\cdots D_{s_1\cdots s_n} .
$$
\hivat{hivx}

\proof
The first statement follows by induction using Definition \cite{hivd}\ (ii)
and the fact that
$$
D_{s_n^{-1}} D_{s_n^{-1}s_{n-1}^{-1} \cdots s_1^{-1} }
\supset
D_{s_n^{-1}} D_{s_n^{-1}s_{n-1}^{-1}} \cdots
  D_{s_n^{-1} \cdots s_1^{-1} } .
$$
By the first statement we have
$$u_{s_n^{-1}\cdots s_1^{-1}}\pi (a^{*})=u_{s_n^{-1}}\cdots u_{s_
1^{-1}}\pi (a^{*}),$$
and the second statement follows from this by taking conjugates.
\eop

\cor
If $(\pi,u,H)$ is a covariant representation, then
$$S=\{u_{s_1}\cdots u_{s_n}:s_1,\ldots,s_n\in G\}$$
is a unital inverse semigroup of partial isometries of H.
\hivat{hivq}

\uj The situation is more delicate than it may appear at
first glance:  the following example shows that a set of
partial isometries with commuting initial and final
projections does not necessarily generate an inverse
semigroup of partial isometries.

\example Let $\alpha$$\in (0,{{\pi}\over 2})$ and
$$U=\left(\matrix{\cos\alpha&-\sin\alpha&0\cr
\sin\alpha&\cos\alpha&0\cr
0&0&0\cr}
\right),\quad V=\left(\matrix{0&0&0\cr
0&\cos\alpha&-\sin\alpha\cr
0&\sin\alpha&\cos\alpha\cr}
\right)$$
be partial isometries on ${\bf C}^3$. A short calculation shows
that $U^2,V^2,UV,VU$ are all partial isometries and so all
the initial and final projections of $U$ and $V$ commute.
But
$$W=(UV)^2=\left(\matrix{0&-\sin\alpha\cos^3\alpha&\sin^2\alpha\cos^
2\alpha\cr
0&\cos^4\alpha&-\sin\alpha\cos^3\alpha\cr
0&0&0\cr}
\right)$$
is not a partial isometry because
$$WW^{*}W=\left(\matrix{0&-\sin\alpha\cos^7\alpha&\sin^2\alpha\cos^
6\alpha\cr
0&\cos^8\alpha&-\sin\alpha\cos^7\alpha\cr
0&0&0\cr}
\right)\ne W.$$
This example is a modification of an idea of Marcelo
Laca.

\newsection{Action of an inverse semigroup}

\uj
In this section we define an action of a unital inverse semigroup and
a covariant representation of such an action. The assumption of the
identity of the semigroup is for technical reasons. In the absence of an
identity
we can easily add one. There is a connection between crossed products by
partial actions and crossed products by semigroup actions, which we are
going to explore in Section~6.

\definition
Let $A$ be a $C^*$-algebra and $S$ be a unital inverse semigroup with identity
$e$.
An {\it action}
of $S$ on $A$ is a semigroup homomorphism
$s\mapsto (\b_s,E_{s^*},E_s):S\to \PAut(A)$, with
$E_e=A$.

\uj
Notice that $\b_{s^*}=\b_s^{-1}$ for all $s\in S$
so the notations $E_{s^*}$ and $E_s$ make sense.
It can be shown as in Proposition \cite{hiva} that $\b_e$ is the identity map
$\iota$ on $A$. Also if $f\in S$ is an idempotent then so is $\b_f$, which
means $\b_f$ is the identity map on $E_{f^*}=E_f$.

\lemma
If \b is an action of the unital inverse semigroup $S$ on $A$ then
$\b_t(E_{t^*}E_{s}) = E_{ts}$ for all $s,t \in S$.
\hivat{hivvn}

\proof The proof follows from the following calculation:
$$\b_t(E_{t^*}E_s)=\b_{t^*}^{-1}(E_{t^*}E_s)=\b_{t^*}^{-1}(E_s)=\hbox{domain }
\b_{s^*}\b_{t^*}=\hbox{domain }\b_{s^*t^*}=E_{ts}.$$
\eop

\example
We have seen in Remark \cite{hivm}\ that a partial action \a of a group
$G$ on $A$ generates a unital inverse subsemigroup
$S=\{
\a_{s_1} \cdots \a_{s_n} : s_1,\ldots, s_n \in G, n\in {\bf Z}
 \}$
of the semigroup
$\PAut(A)$ of partial automorphisms on $A$, and hence determines an action
$s\mapsto\b_s=s$ of the inverse semigroup $S$ on $A$.
\hivat{hivo}

\uj There is a much more important inverse semigroup
action associated with $\alpha$, which we will define in
Section~6 and is based upon the construction used in the
following proposition.

\prop
Let \a be a partial action of a group $G$ on the $C^*$-algebra $A$, and
let $(\pi,u,H)$ be a covariant representation of \a.
Let $S=\{(\a_{g_1}\cdots \a_{g_n},u_{g_1}\cdots u_{g_n}):g_1,\ldots,g_n\in
G\}$.
Then $S$ is a unital inverse semigroup with coordinatewise multiplication.
For $s=(\a_{g_1}\cdots \a_{g_n},u_{g_1}\cdots u_{g_n})\in S$ let
$$
\eqalign{
E_{s^*}&=D_{g_n^{-1}}D_{g_n^{-1}g_{n-1}^{-1}} \cdots D_{g_n^{-1}
  \cdots g_1^{-1}}\cr
E_s&=D_{g_1}D_{g_1g_2} \cdots D_{g_1 \cdots g_n } \cr
\beta_s&=\a_{g_1}\cdots \a_{g_n}:E_{s^*} \to E_s. \cr
}$$
Then \b is an action of $S$ on $A$.
\hivat{hivr}

\proof
By Remark \cite{hivm} and Corollary \cite{hivq}, $S$ is a
unital inverse semigroup with identity $(\a_e,u_e)$,
where $e$ is the identity of $G$.
It is clear that \b is a semigroup homomorphism with
$E_{(\a_e,u_e)}=D_e=A$.
\eop

\definition
Let \b be an action of the unital inverse semigroup
$S$ on $A$. A {\it covariant representation} of \b is
a triple $(\pi,v,H)$ where $\pi:A\to B(H)$ is a nondegenerate representation
of $A$ on the Hilbert space $H$ and $s\mapsto v_s$ is a
semigroup homomorphism from $S$ into an inverse semigroup of
partial isometries on $H$  such that
\nobreak
\items{(i)}$v_s\pi(a)v_{s^{*}}=\pi(\b_s(a))\quad\hbox{for all}\quad
       a\in E_{s^*}$
\items{(ii)} $v_s$ has initial space $\pi(E_{s^*})H$
and final space $\pi(E_s)H$.
\hivat{hivp}

\uj
It can be shown similarly as in the group case (Proposition \cite{hivn})
that $v_e=1_H$ and $v_{s^*}=v_s^*$.
We denote the class of all covariant representations of $(A,S,\b {
)}$ by
$\CovRep(A,S,\b)$.

\prop
Keeping the notations of Proposition \cite{hivr}
define $v:S\to B(H)$ by $v_s=u_{g_1}\cdots u_{g_n}$, where
$s=(\a_{g_1}\cdots \a_{g_n},u_{g_1}\cdots u_{g_n})$.
Then $(\pi,v,H)\in \CovRep(A,S,\b)$.
Furthermore if $(\rho,z,K)\in \CovRep(A,S,\b)$ then the function
$$ w:G\to B(K)\quad \hbox{defined by}\quad w_g=z(\a_g,u_g)
$$
gives a covariant representation $(\rho,w,K)$ of $(A,G,\a)$.
The connections can be visualized by the following diagram
$$

\def\mapright#1{\smash{
    \mathop{\longrightarrow}\limits^{#1}}}
\def\mapleft#1{\smash{
    \mathop{\longleftarrow}\limits^{#1}}}
\def\mapdown#1{\big\downarrow
    \rlap{$\vcenter{\hbox{$\scriptstyle#1$}}$}
    }
\def\mapse#1{_{#1}\searrow}
\def\mapne#1{^{#1}\nearrow}
\def\mapsw#1{\swarrow _{#1}}
\def\mapnw#1{\nwarrow ^{#1}}
\matrix{
               & G  &   \cr
   & \mapsw{u} \quad \mapdown{} \quad \mapse{w}\cr
& B(H) \quad  \mapleft{v}\ \ S \ \ \mapright{z}\quad B(K)  &    \cr
   & \mapnw{\pi} \quad \         \quad \mapne{\rho}\cr
               & A  &   \cr
}
$$

\hivat{hivw}

\def\Ad{\hbox{\rm Ad }}

\proof
It is clear that $v$ is a semigroup homomorphism from $S$ into an inverse
semigroup of partial isometries on $H$. To check Definition
\cite{hivp} (i) let
$s=(\a_{g_1}\cdots \a_{g_n},u_{g_1}\cdots u_{g_n})\in S $ and
$a\in E_{s^*}=
D_{g_n^{-1}}D_{g_n^{-1}g_{n-1}^{-1}} \cdots D_{g_n^{-1} \cdots g_1^{-1}}.$
Using Definition \cite{hivd} (i) and Lemma \cite{hivf}\
we have
$$\eqalign{v_s\pi (a)v_{s^{*}}&=\Ad u_{g_1}\cdots u_{g_n}\circ\pi
(a)\cr
&=\Ad u_{g_1}\cdots u_{g_{n-1}}\circ\pi\circ\a_{g_n}(a)\cr
&=\cdots\cr
&=\pi\circ\a_{g_1}\cdots\a_{g_n}(a)\cr
&=\pi (\b{_s}(a)).\cr}
$$
By Theorem \cite{hiv1} $v_s$ has the desired initial and final spaces.
For the second part of the theorem let $a\in D_{g^{-1}}=E_{s^*}$,
where $s=(\a_g,u_g)$. Then
$$
w_g\rho(a)w_{g^{-1}} = z_s\rho(a)z_{s^*}=\rho(\b_s(a)) =\rho(\a_g(a)),
$$
and so $w$ satisfies Definition \cite{hivd} (i).
To check Definition \cite{hivd} (ii) let $g_1,g_2\in G$,
$h\in \rho(D_{g_2^{-1}}D_{g_2^{-1}g_1^{-1}})K$ and let
$s=(\a_{g_1g_2},u_{g_1g_2}), s_1=(\a_{g_1},u_{g_1}),
s_2=(\a_{g_2},u_{g_2})\in S$.
By Definition \cite{hiv77} (ii)
$\a_{g_1g_2}(\a_{g_1}\a_{g_2})^*=\a_{g_1}\a_{g_2}(\a_{g_1}\a_{g_2})^*$.
By Definition \cite{hivd} (ii) and Theorem \cite{hiv1}
$u_{g_1g_2}(u_{g_1}u_{g_2})^*=u_{g_1}u_{g_2}(u_{g_1}u_{g_2})^*$.
Hence $s(s_1s_2)^*=s_1s_2(s_1s_2)^*$ and so
$z_s z_{(s_1s_2)^*}=z_{s_1s_2} z_{(s_1s_2)^*}$.
Since the final space of $z_{(s_1s_2)^*}$ is
$\rho(D_{g_2^{-1}}D_{g_2^{-1}g_1^{-1}})K$, it follows that
$z_sh=z_{s_1s_2}h$.
Thus
$$w_{g_1g_2}h=z_sh=z_{s_1s_2}h=z_{s_1}z_{s_2}h=w_{g_1}w_{g_2}h$$
as desired. It is clear that $w_g$ has the required initial and final spaces.
\eop

\uj
Notice that if in the previous theorem we let $z=v$ then the construction
gives $w=u$.

\uj
Not every unital inverse semigroup action arises from a partial
action via the construction of Proposition \cite{hivr}, as we can see in
the next example.

\example
Let $S=\{e,f\}$ be the unital inverse semigroup that contains the identity $e$
and an idempotent $f\neq e$.
Let $A={\bf C}$ and $\b_s$ be the identity map $\iota$ of $A$ for all $s\in S$.
Suppose there is a partial action $(A,G,\a)$ and a covariant
representation $(\pi,u,H)$ of \a so that $S$ can be identified with the
inverse semigroup $\{(\a_{g_1}\cdots \a_{g_n},u_{g_1}\cdots u_{g_n}):g_1,
\ldots,g_n\in G\}$ and $\b_s=\a_{g_1}\cdots \a_{g_n}$ for all $s=
(\a_{g_1}\cdots \a_{g_n},u_{g_1}\cdots u_{g_n}) \in S$.
Clearly $e$ is identified with $(\iota,1_H)$, where $1_H$ is the identity of
$B(H)$.
Suppose $f$ is identified with
$(\a_{g_1}\cdots \a_{g_n},u_{g_1}\cdots u_{g_n})$.
By the definition of $\b_s$,
for all $g_1,\ldots,g_n\in G$ we have $\a_{g_1}\cdots \a_{g_n}=\iota$.
Since $f$ is an idempotent $u_{g_1}\cdots u_{g_n}$ is an idempotent too.
Hence for all $h\in H$ we have
$$\eqalign{
h &=\pi(1)(h)=\pi(\b_f(1))(h)\cr
&=(u_{g_1}\cdots u_{g_n}\pi(1)u_{g_1}\cdots u_{g_n})(h)\cr
&=u_{g_1}\cdots u_{g_n}(h).
}$$
This means that $u_{g_1}\cdots u_{g_n}$ must be the identity
of $B(H)$. But this is a contradiction since $e$ and $f$ are different
elements of $S$.
\hivat{hiv4}

\newsection{The crossed product} \uj
McClanahan[\cite{mcc}] defines the {\it partial crossed }
{\it product} $A\times_{\a}G$ of the $C^{*}$-algebra $A$ and the group $
G$ by
the partial action $\a$ as the enveloping $C^{*}$-algebra of
$L=\{x\in l^1(G,A):x(g)\in D_g\}$ with multiplication and
involution
$$\eqalign{(x*y)(g)&=\sum\limits_{h\in G}\a_h[\a_{h^{-1}}(x(h))y(
h^{-1}g)],\cr
x^{*}(g)&=\a_g(x(g^{-1})^{*}).\cr}
$$

He shows that there is bijective
correspondence $(\pi,u,H)\leftrightarrow (\pi\times u,H)$ between
covariant representations of $(A,G,\a)$ and nondegenerate representations
of $A\times_\a G$, where $\pi\times u$ is the extension of the
representation of $L$ defined by
$$x\mapsto\sum_{g\in G}\pi(x(g))u_g .$$
We are going to follow his footsteps constructing
the {\it crossed product} of a $C^*$-algebra and a unital inverse semigroup by
an action \b.

Let $\b$ be an action of the unital inverse semigroup $S$ on the
$C^{*}$-algebra $A$. Let
$$L=\{x\in l^1(S,A):x(s)\in E_s\}$$
have the norm, scalar multiplication and addition inherited from
$l^1(S,A)$. Define multiplication and involution on $L$ by
$$\eqalign{
(x*y)(s) &= \sum\limits_{rt=s}\b_r[\b_{r^*}(x(r))y(t)] \cr
x^*(s)   &= \b_s(x(s^*)^*) .\cr
}$$
Notice that by Lemma \cite{hivvn} $(x*y)(s)\in E_s$.
A short calculation shows that $\|x*y\|\leq \|x\|\|y\|$ and so
$x*y \in L$ :
$$\eqalign{\|x*y\|&=\sum_{s\in S}\sum_{rt=s}\|\b{_r}(\b{_{r^{*}}}
(x(r))y(t))\|\cr
&=\sum_{s\in S}\sum_{rt=s}\|\b{_{r^{*}}}(x(r))y(t)\|\cr
&\leq\sum_{r\in S}\sum_{t\in S}\|\b{_{r^{*}}}(x(r))\|\|y(t)\|\cr
&\leq\sum_{r\in S}\|\b{_{r^{*}}}(x(r))\|\sum_{t\in S}\|y(t)\|\cr
&\leq\|x\|\|y\|.\cr}
$$
One can easily check that $\|x^{*}\|=\|x\|$ and so $x^{*}\in L$.  We
are going to denote by $a\delta_s$ the function in $L$ taking the
value $a$ at $s$ and zero at every other element of $S$.
Notice that $a_s\delta_s*a_t\delta_t=\b{_s}(\b{_{s^{*}}}(a_s)a_t)
\delta_{st}$ and
$(a\delta_s)^{*}=\beta_{s^{*}}(a^{*})\delta_{s^{*}}$.

\prop
$L$ is a Banach *-algebra.

\proof
Let $x,y,z\in L$ and $a\in{\bf C}$. Then it is easy to see that
$(x+y)^*=x^*+y^*$ and $(ax)^*=\bar ax^*$. We also have $x^{**}=x$ since
$$x^{**}(s)=\b{_s}(x^{*}(s^{*})^{*})=\b{_s}(\b{_{s^{*}}}(x(s^{**}
)^{**}))=x(s).$$
Next we show that $(x*y)^*=y^**x^*$. It suffices to show this for
$x=a_s\d_s$ and $y=a_t\d_t$. We have
$$\eqalign{
(a_s\d_s*a_t\d_t)^*
&= (\b_s(\b_{s^*}(a_s)a_t)\d_{st})^* \cr
&= \b_{t^*s^*}((\b_s(\b_{s^*}(a_s)a_t))^*)\d_{t^*s^*} \cr
&= \b_{t^*s^*}(\b_s(a_t^*\b_{s^*}(a_s^*)))\d_{t^*s^*} \cr
&= \b_{t^*}(a_t^*\b_{s^*}(a_s^*))\d_{t^*s^*} \cr
&= \b_{t^*}(\b_t(\b_{t^*}(a_t^*))\b_{s^*}(a_s^*))\d_{t^*s^*} \cr
&= \b_{t^*}(a_t^*)\d_{t^*}*\b_{s^*}(a_s^*)\d_{s^*} \cr
&= (a_t\d_{t})^**(a_s\d_{s})^* \cr
}$$
as desired.

Finally we show that $(x*y)*z=x*(y*z)$.  Again it
suffices to show this for $x=a_r\d{_r}$, $y=a_s\d{_s}$ and $z=a_t
\d{_t}$.
If $\{u_{\lambda}\}$ is an approximate identity for $E_{s^{*}}$, then we
have
$$\eqalign{(a_r\d{_r}*a_s\d{_s})*a_t\d{_t}&=\b{_r}(\b{_{r^{*}}}(a_
r)a_s)\d{_{rs}}*a_t\d{_t}\cr
&=\b{_{rs}}(\b{_{s^{*}r^{*}}}(\b{_r}(\b{_{r^{*}}}(a_r)a_s))a_t)\d{_{
rst}}\cr
&=\lim_{\lambda}\b{_{rs}}(\b{_{s^{*}}}(\b{_{r^{*}}}(a_r)a_s)u_{\lambda}
a_t)\d{_{rst}}\cr
&=\lim_{\lambda}\b{_r}(\b{_{r^{*}}}(a_r)a_s\b{_s}(u_{\lambda}a_t)
)\d{_{rst}}\cr
&=\lim_{\lambda}\b{_r}(\b{_{r^{*}}}(a_r)\b{_s}(\b{_{s^{*}}}(a_s)u_{
\lambda}a_t))\d{_{rst}}\cr
&=\b{_r}(\b{_{r^{*}}}(a_r)\b{_s}(\b{_{s^{*}}}(a_s)a_t))\d{_{rst}}\cr
&=a_r\d{_r}*\b{_s}(\b{_{s^{*}}}(a_s)a_t)\d{_{st}}\cr
&=a_r\d{_r}*(a_s\d{_s}*a_t\d{_t})\cr}
$$
as desired.
\eop

\definition
If $(\pi,v,H)\in \CovRep(A,S,\b)$ then define $\pi\times v:L\to B(H)$ by
$$ (\pi\times v)(x)=\sum\limits_{s\in S} \pi(x(s))v_s.
$$

\prop
$(\pi\times v)$ is a *-homomorphism onto the $C^*$-algebra
$$C^*(\pi,v)=\overline{\sum_{s\in S}\pi(E_s)v_s} . $$
\hivat{lasthiv}

\proof First notice that by Definition \cite{hivp},
$v_sv_{s^{*}}\pi (a_s)=\pi (a_s)=\pi (a_s)v_sv_{s^{*}}$ for all $
a_s\in E_s$.  If $s,t\in S$,
$a\in E_s$ and $b\in E_t$ then
$$\eqalign{\pi (a)v_s\pi (b)v_t&=v_s\pi (\b{_{s^{*}}}(a))v_{s^{*}}
v_s\pi (b)v_t\cr
&=v_s\pi (\b{_{s^{*}}}(a)b)v_t\cr
&=v_s\pi (\b{_{s^{*}}}\b{_s}(\b{_{s^{*}}}(a)b))v_t\cr
&=v_sv_{s^{*}}\pi (\b{_s}(\b{_{s^{*}}}(a)b))v_sv_t\cr
&=\pi (\b{_s}(\b{_{s^{*}}}(a)b))v_{st}\cr}
$$
and
$$(\pi (a)v_s)^{*}=v_{s^{*}}\pi (a^{*})=v_{s^{*}}v_s\pi (\b{_{s^{
*}}}(a^{*}))v_{s^{*}}=\pi (\b{_{s^{*}}}(a^{*}))v_{s^{*}},$$
which shows that $C^*(\pi,v)$ is really a $C^*$-algebra.
It is clear that $\pi\times v$ is linear.
It suffices to verify the multiplicativity of $\pi\times v$ for elements
of the form $a_s\delta_s$. We have
$$\eqalign{
(\pi\times v)(a_s\delta_s*a_t\delta_t)
 &=(\pi\times v)(\b_s(\b_{s^*}(a_s)a_t)\delta_{st}) \cr
 &=\pi(\b_s(\b_{s^*}(a_s)a_t))v_{st} \cr
 &=\pi(a_s)v_s\pi(a_t)v_t  \cr
 &=(\pi\times v)(a_s\delta_s)(\pi\times v)(a_t\delta_t)
}$$
as desired.
The following calculation shows that $\pi\times v$ preserves the *-operation:
$$\eqalign{
(\pi\times v)(a_s\delta_s)^*
 &=(\pi(a_s)v_s)^* \cr
 &=\pi(\b_{s^*}(a_s^*)) v_{s^*} \cr
 &=(\pi\times v)(\b_{s^*}(a_s^*)\delta_{s^*}) \cr
 &=(\pi\times v)((a_s\delta_s)^*) . \cr
}$$
It is clear that $(\pi\times v)$ is onto.
\eop

\definition
\hivat{crospdef}
Let $A$ be a $C^*$-algebra and \b be an action of the unital inverse semigroup
$S$ on $A$.
Define a seminorm $\|.\|_c$ on $L$ by
$$
\|x\|_c=\sup\{\|(\pi\times v)(x)\|:(\pi,v)\in \CovRep(A,S,\b) \} .
$$
Let $I=\{x\in L:\|x\|_c=0 \}$.
The \it crossed product \rm $A\times_\b S$ of the $C^*$-algebra $A$ and the
semigroup $S$ by the action \b is the $C^*$-algebra gotten by the completion of
the quotient $L/I$ with  respect to $\|.\|_c$. We denote the quotient map by
$\Phi$.

\def\l {\lambda}
\def\d {\delta}

\uj
Since $\Phi(L)$ is dense in $A\times_\b S$,
it is clear that $\pi\times v$ induces a nondegenerate
representation of $A\times_\b S$.
We denote this representation also by $\pi\times v$.

The following lemma shows that the ideal $I$ may be nontrivial:

\lemma If $s,t\in S$ so that $s\leq t$, that is, $s=ft$ for some
idempotent $f\in S$ then $\Phi (a\delta_s)=\Phi (a\delta_t)$ for all $
a\in E_s$.  In
particular $\Phi (a\delta_s)=\Phi (a\delta_e)$ if $s$ is an idempotent.
\hivat{hiv2}

\proof\ It is clear that $a\in E_t$.  If $(\pi,v)\in\CovRep(A,S,
\beta )$
then
$$(\pi\times v)(a\delta_{ft}-a\delta_t)=\pi (a)v_fv_t-\pi (a)v_t=
0,$$
which shows $\Phi (a\delta_s-a\delta_t)=0.$ The second statement
follows from the fact that $s=se$.  \eop

\uj
If $(\Pi,H)$ is a representation of $A\times_\b S$ and $x\in L$ then we are
going to write $\Pi (x)$ instead of the more precise $\Pi (\Phi (
x))$.

\prop
Let $(\Pi,H)$ be a nondegenerate representation of $A\times_\b S$. Define a
representation $\pi$ of $A$ on $H$ by
$$ \pi(a)= \Pi(a\d_e).
$$
Let $v:S\to B(H)$ defined by
$$ v_s=\slim\limits_\l\Pi(u_\l\d_s)
$$
where $\{u_\l\}$ is an approximate identity for $E_s$ and $\slim_\l$ denotes
the
strong operator limit. Then $(\pi,v,H)\in \CovRep(A,S,\b)$.

\proof
First we show that $v_s$ is well defined.
If $h\in \pi(E_{s^*})H$ then $h=\Pi(a\d_e)k$ for some $a\in E_{s^*}$ and
$k\in H$. So
$$\eqalign{\lim_{\l}\Pi (u_{\l}\d{_s})(h)&=\lim_{\l}\Pi (u_{\l}\d{_
s})\Pi (a\d{_e})(k)\cr
&=\lim_{\l}\Pi (u_{\l}\d{_s}*a\d{_e})(k)\cr
&=\lim_{\l}\Pi (\b{_s}(\b{_{s^{*}}}(u_{\l})a)\d{_s})(k)\cr
&=\Pi (\b{_s}(a)\d{_s})(k),\cr}
$$
where we used the fact that $\b{_{s^{*}}}(u_{\l})$ is an approximate
identity for $E_{s^{*}}$.  Note that the limit is independent of
the choice of $\{u_{\lambda}\}$ since the expression $h=\pi (a\delta_
e)k$ was.
On the other hand if $\langle h,\pi (E_{s^{*}})H\rangle =\langle
h,\Pi (E_{s^{*}}\d{_e})H\rangle =0$ then
$$\eqalign{\lim_{\l}\Pi (u_{\l}\d{_s})(h)&=\lim_{\l}\Pi (\b{_s}(\b{_{
s^{*}}}(\sqrt {u_{\l}})\b{_{s^{*}}}(\sqrt {u_{\l}}))\d{_s})(h)\cr
&=\lim_{\l}\Pi (\sqrt {u_{\l}}\d{_s}*\b{_{s^{*}}}(\sqrt {u_{\l}})
\d{_e})(h)\cr
&=\lim_{\l}\Pi (\sqrt {u_{\l}}\d{_s})\Pi (\b{_{s^{*}}}(\sqrt {u_{
\l}})\d{_e})(h).\cr}
$$
But $\langle \Pi(\b_{s^*}(\sqrt{u_\l})\d_e)(h),H \rangle=
 \langle h,\Pi(\b_{s^*}(\sqrt{u_\l})\d_e) H \rangle=0$ and so
$\lim_\l\Pi(u_\l\d_s)(h)=0$. Hence $v_s$ is well defined.
It is easy to see that $v_s$ is a bounded linear transformation.
The following calculation shows that $v_s^*=v_{s^*}$:
$$ \eqalign{
v_s^* &= \slim_\l\Pi(u_\l\d_s)^*  \cr
      &= \slim_\l\Pi(\b_{s^*}(u_\l)\d_{s^*})
      = v_{s^*} . \cr
}$$
In order to see that $v_s$ is a partial isometry we need to show that
$v_sv_s^*v_s=v_s$. It suffices to
show that $v_s^*v_sh=h$ for $h\in \pi(E_{s^*})H$, since we have seen above
that $(\pi(E_{s^*})H)^{\perp}\subset \hbox{ker\ } v_s$.
Let $h=\Pi(a\d_e)k$ where $a\in E_{s^*}$ and $k\in H$.
Using the fact that $\Phi(a\d_{s^*s})=\Phi(a\d_e)$ we have
$$\eqalign{
 v_s^*v_s(h) &= \lim_\l\Pi(u_\l\d_{s^*})\Pi(\b_s(a)\d_s)(k) \cr
   &= \lim_\l\Pi(\b_{s^*}(\b_s(u_\l)\b_s(a))\d_{s^*s})(k) \cr
   &= \Pi(a\d_{s^*s})(k) = \Pi(a\d_e)(k)=h . \cr
}$$
This calculation also shows that $v_s$ has initial space
$v_s^*v_s(H)=\pi(E_{s^*})H$. A similar calculation shows that $v_s$ has
final space $\pi(E_s)H$.
To see that $v$ is a semigroup homomorphism let $h=\Pi(a\d_e)k\in
\pi(E_{t^*s^*})H$ where
$a\in E_{t^*s^*}$ and $k\in H$. Then
$$\eqalign{
v_sv_t(h) &=v_sv_t\Pi(a\d_e)(k)\cr
  &= v_s\Pi(\b_t(a)\d_t)(k),\qquad \hbox{by Lemma \cite{hivvn}}\cr
  &= \lim_\l \Pi(u_\l\d_s)\Pi(\b_t(a)\d_t)(k) \cr
  &= \lim_\l \Pi(\b_s(\b_{s^*}(u_\l)\b_t(a)\d_{st})(k) \cr
  &= \Pi(\b_s\b_t(a)\d_{st})(k),\qquad \hbox{since } \b_t(a)\in E_{s^*}\cr
  &= v_{st}\Pi(a\d_e)(k) = v_{st}(h) .
}$$
On the other hand if $\langle h,\Pi(E_{t^*s^*}\d_e)H \rangle=0$ then
$v_{st}(h)=0$.
We show $v_sv_t(h)=0$ as well.
If $w_\mu$ is an approximate identity for $E_t$ then
$$\eqalign{
v_sv_t(h) &= \slim_\l\Pi(u_\l\d_s)\slim_\mu\Pi(w_\mu\d_t)(h)  \cr
 &= \slim_{\mu,\l}\Pi(u_\l\d_s*w_\mu\d_t)(h) \cr
 &= \slim_{\mu,\l}\Pi(\b_s(\b_{s^*}(u_\l)w_\mu)\d_{st})(h) .
}$$
By Lemma \cite{hivvn}  $\b{_s}(\b{_{s^{*}}}(u_{\l})w_{\mu})\in E_{
st}$,
and so can be factored as $xy$ with $x,y\in E_{st}$ (by the Cohen-Hewitt
factorization theorem). Now we have
$$\Pi (xy\d{_{st}})(h)=\Pi (x\d{_{st}})\Pi (\b{_{t^{*}s^{*}}}(y)\d{_
e})(h).$$
But $\langle\Pi (\b{_{t^{*}s^{*}}}(y)\d{_e})(h),H\rangle =\langle
h,\Pi (\b{_{t^{*}s^{*}}}(y)\d{_e})H\rangle =0$
and so $v_sv_t(h)=0$.

For the covariance condition, if $a\in E_{s^{*}}$ then
$$\eqalign{
v_s\pi(a)v_{s^*}
  &= \slim_{\l,\mu}\Pi(u_\mu\d_s)\Pi(a\d_e)\Pi(\b_{s^*}(u_\l)\d_{s^*}) \cr
  &= \slim_{\l,\mu}\Pi(u_\mu\b_s(a)u_\l\d_{ss^*}) \cr
  &= \Pi(\b_s(a)\d_{e}) \cr
}$$
as desired. The nondegeneracy of $\pi$ follows from that of $\Pi$,
since $\{u_{\lambda}\delta_e\}$ is an approximate identity for $A
\times_{\b{}}S$
whenever $\{u_\lambda\}$ is an approximate identity for $A$.
\hbox{\ \ \ \ \ \ \ \ \ \ }
\eop

\prop
The correspondence $(\pi,v,H)\leftrightarrow (\pi\times v,H)$ is a bijection
between $\CovRep(A,S,\b)$ and the class of nondegenerate representations of
$A\times_\b S$.
\hivat{hivz}

\proof
\def\t{\tilde}
Let $(\t\pi,\t v,H)\in \CovRep(A,S,\b)$. Let $(\pi,v)$ the the covariant
representation induced by $\t\pi\times\t v$. Then
$$ \pi(a)=(\t\pi\times\t v)(a\d_e)=\t\pi(a)v_e=\t\pi(a)
$$
$$ v_s=\slim_\l(\t\pi\times\t v)(u_\l\d_s)=\slim_\l\t\pi(u_\l)\t v_s=\t v_s .
$$
The last equality holds because $\t\pi(u_\l)$ converges strongly to the
identity on $\t\pi(E_s)H$.
On the other hand if $\Pi$ is a representation of $A\times_\b S$,
$(\pi ,v)$ is induced by $\Pi$ and $a\in E_s$ then
$$\eqalign{(\pi\times v)(a\delta_s)&=\pi (a)v_s\cr
&=\Pi (a\delta_e)\slim_{\lambda}\Pi(u_{\lambda}\delta{_s})\cr
&=\slim_{\l}\Pi (a\d{_e}*u_{\l}\d{_s})\cr
&=\slim_{\l}\Pi (au_{\l}\d{_s})\cr
&=\Pi (a\d{_s}).\cr}
$$
Thus the correspondence is a bijection.
\eop

\prop If $\beta$ is an action of the semilattice $S$ on a $C^{*}$-algebra
$A$. Then $A\times_{\beta}S$ is isomorphic to $A$.
\hivat{hivlatt}

\proof If $(\pi ,v)$ is any covariant representation of $\beta$ then
by Lemma \cite{hiv2} ($\pi\times v)(A\times_{\beta}S)=\pi (A)$.  In particular
if $\pi$ is a faithful representation of $A$ on the Hilbert
space $H$ and $v_f$ is the projection onto $\pi (E_f)H$ for all
$f\in S$ then $(\pi ,v)$ is a covariant representation of $\beta$ and
$(\pi\times v)(A\times_{\beta}S)=\pi (A)\cong A$.  Every representation of $
A\times_{\beta}S$
factors through $(\pi\times v)$ for if $(\rho\times z)$ is a representation of
$A\times_{\beta}S$ then $(\rho\times z)=\rho\circ\pi^{-1}\circ (\pi
\times v)$.  Thus $(\pi\times v)$ is
faithful and so $A\times_{\beta}S$ and $A$ must be isomorphic.  \eop

\uj
The following example shows that unlike in the partial action case the
crossed product $A\times_\b S$ is
not the enveloping $C^{*}$-algebra of $L$ in general.

\example
Let $S=\{e,f\}$ be the unital inverse semigroup that contains the identity $e$
and an idempotent $f$.
Let $A={\bf C}$ and $\b_s$ be the identity map of $A$ for all $s\in S$ as in
Example \cite{hiv4}.
It is clear that $L=l^1(S)$. Wordingham [\cite{wordi}] shows that the
left regular representation
of $l^1(S)$ on $l^2(S)$ is faithful and so the enveloping $C^*$-algebra
cannot be the same as $A\times_\b S$, which is isomorphic to
{\bf C} by Proposition \cite{hivlatt}.
\hivat{hiv3}

\def \d {\delta}

\def \P {\Psi}
\def \CC {\hbox{\bf C}}

\uj The next two results describe two quite different
crossed products associated with an inverse semigroup
itself.

\prop
Let $S$ be a unital inverse semigroup, and let $\b_s$ be the identity map of
{\bf C} for all $s\in S$. Then $\b_s$ is an action of $S$ on \CC, and
$\hbox{\bf C}\times_\b S$ is isomorphic to the $C^*$-algebra of the maximal
group homomorphic image of $S$.

\proof
For $s,t\in S$ let $s\sim t$ if and only if there is an idempotent $
f\in S$
so that $fs=ft$. Then $\sim$ is a congruence on $S$ and $G=S/\sim$ is the
maximal group homomorphic image of $S$. Let $[s]$ denote
the equivalence class of $s\in S$.

Let $\Phi:L\to L/I$ be as in Definition \cite{crospdef}.
Define a *-homomorphism
$$\P :G\to L/I\quad\hbox{\rm by}\quad\P ([s])=\Phi (\delta_s).$$
$\P$ is well defined since if $s\sim t$ then there is an
idempotent $f\in S$ so that $fs=ft$ and so by Lemma
\cite{hiv2} we have $\Phi (\d{_s})=\Phi (\d{_f})\Phi (\d{_s})=\Phi
(\d{_f})\Phi (\d{_t})=\Phi (\d{_t}).$ $\P$
extends to a *-homomorphism $\P :C^{*}(G)\to\hbox{\rm {\bf C}}\times_{
\beta}S$.  The
extension is onto since $\P (G)=\Phi (\delta_S)$ has dense span in
$\hbox{\rm {\bf C}}\times_{\beta}S$.

Going the other way, define a covariant representation of
$(\CC,S,\beta )$ in $C^{*}(G)$, considered to be represented on a
Hilbert space, by
$$\eqalign{\pi :\hbox{\rm {\bf C}}\to C^{*}(G),\quad&\pi (a)=ae\cr
v:S\to C^{*}(G),\quad&v_s=[s],\cr}
$$
where $G$ is identified with its canonical image in $C^{*}(G)$.
In fact, $(\pi ,v)$ satisfies the covariant condition
$v_s\pi (a)v_{s^{*}}=a[ss^{*}]=\pi (a)$
because all the idempotents are congruent.

The covariant representation $(\pi,v)$ defines a *-homomorphism
$$(\pi\times v):\CC\times_{\b{}}S\to C^{*}(G).$$
For $[s]\in G$ we have
$$(\pi\times v)\circ\P ([s])=(\pi\times v)(\delta_s)=v_s=[s].$$
So $(\pi\times v)\circ\P$ is the identity map on $G$. Since a representation
of $G$ has a unique extension to $C^{*}(G)$, $(\pi\times v)\circ\P$ must be the
identity map on $C^{*}(G)$ and hence $\P$ is an isomorphism between $
C^{*}(G)$ and
$\CC\times_\b S$.
\eop

\prop
Let $S$ be a unital inverse semigroup with idempotent semilattice
$E$. Define a semigroup action $\b$ of $S$ on $C^*(E)$ so that
$\b_s:E_{s^*}\rightarrow E_s$ is determined by $\b_s(\d_f)=\d_{sfs^*}$,
where $E_s$ is the closed span of the set $\{\d{_f}:f\leq ss^{*}\}$ in $
C^{*}(E)$.
Then $C^*(S)$ is isomorphic to $C^*(E)\times_\beta S$.

\proof Let $\Phi :L\to L/I$ be as in Definition \cite{crospdef}.  We
are going to identify $S$ with its canonical image in
$C^{*}(S)$.  Define a *-homomorphism
$$\P :S\to L/I\quad\hbox{\rm by}\quad\P (s)=\Phi (ss^{*}\delta_s)
.$$
In fact $\P$ is a *-homomorphism since
$$\eqalign{\P (s)\P (t)&=\Phi (ss^{*}\delta_s)\Phi (tt^{*}\delta_
t)=\Phi (\beta_s(\beta_{s^{*}}(ss^{*})tt^{*})\delta_{st})\cr
&=\Phi (\beta_s(s^{*}ss^{*}stt^{*})\delta_{st})=\Phi (ss^{*}stt^{
*}s^{*}\delta_{st})\cr
&=\Phi (st(st)^{*}\delta_{st})=\P (st),\cr}
$$
and
$$\eqalign{\P (s)^{*}&=\Phi ({{ss^{*}}}\delta_s)^{*}=\Phi (\beta_{
s^{*}}(({{ss^{*}}})^{*})\delta_{s^{*}})\cr
&=\Phi (s^{*}ss^{*}s\delta_{s^{*}})=\P (s^{*}).\cr}
$$

$\P$ extends to a *-homomorphism $\P :C^{*}(S)\to C^{*}(E)\times_{
\b{}}S$.  The
extension is onto since if $f\in E$ with $f\leq ss^{*}$ then by
Lemma \cite{hiv2} we have
$$\eqalign{\Phi (f\delta_s)&=\Phi (f\delta_e)\Phi (f\delta_s)=\Phi
(f\delta_f)\Phi (f\delta_s)\cr
&=\Phi (f\delta_{fs})=\Phi (fs(fs)^{*}\delta_{fs})\cr
&=\Psi (fs),\cr}
$$
which means the span of $\P (S)$ is dense in $C^{*}(E)\times_{\beta}
S$.

We show that $\P$ is in fact an isomorphism.  Let
$\pi :C^{*}(E)\to C^{*}(S)$ be the canonical inclusion map and define
$v:S\to C^{*}(S)$ by $v_s=s$.  We show that $(\pi ,v)$ is a covariant
representation of $(C^{*}(E),S,\beta )$ if $C^{*}(S)$ is considered to be
represented on a Hilbert space.  It is clear that $v$ is a
semigroup homomorphism into an inverse semigroup of
partial isometries.  The requirements for the initial and
final spaces are satisfied since $v_sv_s^{*}=ss^{*}=\pi (p_s)$.  It
suffices to verify the covariance condition for elements
$f\in E_{s^{*}}$ where $f\leq s^{*}s$.  For such $f$ we have
$$v_s\pi (f)v_{s^{*}}=sfs^*=\beta{_s}(f),$$
as desired.

The covariant representation $(\pi,v)$ gives a *-homomorphism
$$(\pi\times v):C^*(E)\times_\b S\to C^*(S).$$
Since for $s\in S$ we have
$$(\pi\times v)\circ\P (s)=(\pi\times v)(\Phi (ss^{*}\delta_s))=\pi
(ss^{*})v_s=ss^{*}s=s,$$
the composition $(\pi\times v)\circ\P$ is the identity map on $S$.
Since a representation of $S$ has a unique extension to
$C^{*}(S)$, $(\pi\times v)\circ\P$ must be the identity map on $C^{
*}(S)$, and
so $\P$ is an isomorphism between $C^{*}(S)$ and $C^{*}(E)\times_{
\b{}}S$.
\eop

\newsection{Connection between the crossed products}

\uj
Let $(A,G,\a)$ be a partial action of the group $G$ on $A$. Let $(\Pi,H)$ be a
faithful nondegenerate representation of the crossed product $A\times_\a G$.
McClanahan
[\cite{mcc}] showed that $\Pi=\pi\times u$ for some covariant representation
$(\pi,u,H)$ of $(A,G,\a)$.
Define, as in Proposition \cite{hivr}, a
unital inverse semigroup
$$S=\{(\alpha_{g_1}\cdots\alpha_{g_n},u_{g_1}\cdots u_{g_n}):g_1,
\ldots ,g_n\in G\}$$
and an action $\beta$ of $S$ so that $\beta_s=\alpha_{g_1}\cdots\alpha_{
g_n}$ for
$s=(\alpha_{g_1}\cdots\alpha_{g_n},u_{g_1}\cdots u_{g_n})\in S$.  We are going
to show that
the crossed products $A\times_{\a}G$ and $A\times_{\b{}}S$ are canonically
isomorphic.  First we need the following

\lemma Let $(\rho ,z,K)\in\CovRep(A,S,\b {)}$ and define a covariant
representation $(\rho ,w,K)$ of $(A,G,\a)$ by $w_g=z(\alpha_g,u_g
)$ as in
Proposition \cite{hivw}.  Then
$(\rho\times z)(A\times_{\b{}}S)=(\rho\times w)(A\times_{\a}G)$.  \hivat{hivy}

\proof
By Proposition \cite{lasthiv}, it suffices to show that
$$\sum_{s\in S}\rho (E_s)z_s=\sum_{g\in G}\rho (D_g)w_g.$$
Let $g\in G$ and $s=(\alpha_g,u_g)$.  If $a\in D_g=E_s$ then
$\rho (a)w_g=\rho (a)z_s$ and so $\sum_{s\in S}\rho (E_s)z_s\supset
\sum_{g\in G}\rho (D_g)w_g$.  On
the other hand if $s=(\alpha_{g_1}\cdots\alpha_{g_n},u_{g_1}\cdots
u_{g_n})$ and
$a\in E_s=D_{g_1}D_{g_1g_2}\cdots D_{g_1\cdots g_n}$ then by Corollary
\cite{hivx}
we have
$$\eqalign{\rho (a)z_s&=\rho (a)z(\alpha_{g_1},u_{g_1})\cdots z(\alpha_{
g_1},u_{g_n})\cr
&=\rho (a)w_{g_1}\cdots w_{g_n}=\rho (a)w_{g_1\cdots
g_n}\cr}
$$
and so $\sum_{s\in S}\rho (E_s)z_s\subset\sum_{g\in G}\rho (D_g)w_
g$.
\eop

\uj
Now we can prove our main result.

\theorem\hivat{hi1} Let \a be a partial action of a group $G$
on the $C^{*}$-algebra $A$, and let $(\pi ,u,H)$ be a covariant
representation of $(A,G,\a)$ so that $\pi\times u$ is a faithful
representation of $A\times_{\a}G$.  Define an inverse semigroup
by
$$S=\{(\alpha_{g_1}\cdots\alpha_{g_n},u_{g_1}\cdots u_{g_n}):g_1,
\ldots ,g_n\in G\}$$
and an action $\beta$ of $S$ by $\beta_s=\alpha_{g_1}\cdots\alpha_{
g_n}$ for
$s=(\alpha_{g_1}\cdots\alpha_{g_n},u_{g_1}\cdots u_{g_n})$.  Then the crossed
products
$A\times_{\alpha}G$ and $A\times_{\beta}S$ are isomorphic.

\proof
Let $v_s=u_{g_1}\cdots u_{g_n}$ for $s=(\alpha_{g_1}\cdots\alpha_{
g_n},u_{g_1}\cdots u_{g_n})$.
We know from Proposition \cite{hivw} that $(\pi,v,H)\in \CovRep(A,S,\b)$.
Note that $C^*(\pi,u)=C^*(\pi,v)$.
It suffices to show that $\pi\times v$ is a faithful representation of
$A\times_\b S$ because then by Lemma \cite{hivy}
$(\pi\times v)^{-1}\circ(\pi\times u) : A\times_\a G\to A\times_\b S$ is
an isomorphism. To show this consider another representation of
$A\times_\b S$, which by Proposition \cite{hivz} must be in the form
$\rho\times z$ for some $(\rho,z)\in \CovRep(A,S,\b)$.
By Proposition \cite{hivw} the definition
$w_g=z(\alpha_g,u_g)$ gives a covariant representation $(\rho ,w,
K)$ of $(A,G,\a)$.
By Lemma \cite{hivy} $C^{*}(\rho ,w)=C^{*}(\rho ,z)$. Since $\pi\times
u$ is a
faithful representation, there is a homomorphism $\Theta$ making the following
diagram commute:
$$

\def\mapright#1{\smash{
    \mathop{\longrightarrow}\limits^{#1}}}
\def\mapdown#1{\Big\downarrow
    \rlap{$\vcenter{\hbox{$\scriptstyle#1$}}$}
    }
\def\mapse#1{_{#1}\searrow}
\def\mapsw#1{\swarrow _{#1}}
\matrix{
               & A\times_\a G  &   \cr
   & \mapsw{\pi\times u} \qquad  \mapse{\rho\times w  }  \cr
& C^*(\pi,u) \qquad   \mapright{\Theta} \qquad C^*(\rho,w)  &    \cr
}
$$

We are going to show that the diagram
$$

\def\mapright#1{\smash{
    \mathop{\longrightarrow}\limits^{#1}}}
\def\mapdown#1{\Big\downarrow
    \rlap{$\vcenter{\hbox{$\scriptstyle#1$}}$}
    }
\def\mapse#1{_{#1}\searrow}
\def\mapsw#1{\swarrow _{#1}}
\matrix{
               & A\times_\b S  &   \cr
   & \mapsw{\pi\times v} \qquad  \mapse{\rho\times z  }  \cr
& C^*(\pi,v) \qquad   \mapright{\Theta} \qquad C^*(\rho,z)  &    \cr
}$$
commutes, which is going to prove that $\pi\times v$ is faithful.
Since finite sums of
elements of the form $a\d_s\in L$ (more precisely the images under $\Phi$)
are dense in $A\times_{\b{}}S$, it suffices to check that
$$
\Theta\circ (\pi\times v)(a\d_s)=(\rho\times z)(a\d_s)
$$
where
$s=(\alpha_{g_1}\cdots\alpha_{g_n},u_{g_1}\cdots u_{g_n})$ and $a
\in E_s=D_{g_1}D_{g_1g_2}\cdots D_{g_1\cdots g_n}$.
But this is true by the following calculation:
$$\eqalign{\Theta ((\pi\times v)(a\d{_s}))&=\Theta (\pi (a)v_s)\cr
&=(\rho\times w)\circ (\pi\times u)^{-1}(\pi (a)u_{g_1}\cdots u_{
g_n})\cr
&=(\rho\times w)\circ (\pi\times u)^{-1}(\pi (a)u_{g_1\cdots g_n}
)\cr
&=(\rho\times w)(a\delta_{g_1\cdots g_n})\cr
&=\rho (a)w_{g_1\cdots g_n}\cr
&=\rho (a)w_{g_1}\cdots w_{g_n}\cr
&=\rho (a)z(\alpha_{g_1},u_{g_1})\cdots z(\alpha_{g_n},u_{g_n})\cr
&=\rho (a)z_s\cr
&=(\rho\times z)(a\delta_s),\cr}
$$
where we have appealed to Corollary \cite{hivx} twice more.
\eop

\def\C{{\bf C}^2}

\example Let $A=\C$, $G={\bf Z}$, $D_0=A$, $D_{-1}=\{(a,0):a\in A
\}$,
$D_1=\{(0,a):a\in A\}$ and $D_n=\{(0,0)\}$ for $n\in G\setminus \{
-1,0,1\}$.  Let
$\alpha_0$ be the identity map $\a_1$ be the forward shift
$\a_1(a,0)=(0,a)$ and define $\a_n=\a_1^n$ for $n\ne 0$.  Then
$A\times_{\a}G$ is isomorphic to the matrix algebra $M_2$
[\cite{ex}], [\cite{mcc}, Example 2.5].

We construct a faithful
representation $\pi\times u$ of the partial crossed product
\hbox{$A\times_\a G$}.
Let $\pi$ be the representation of $A$ on the Hilbert space $H=\C$ by
multiplication operators, that is,
$$
\pi(a_1,a_2)(h_1,h_2)=(a_1h_1,a_2h_2)
$$
for $(a_1,a_2)\in A$ and $(h_1,h_2)\in H$. Let $u_1$ be the forward shift,
$u_{-1}$ the backward shift on $H$, and let $u_n$ be the constant zero map for
all $n\in G\setminus\{-1,0,1\}$.

The unital inverse semigroup $S$ generated by
$\{(\alpha_n,u_n):n\in G\}$ contains six elements
$$S=\{e,f,s,s^{*},s^{*}s,ss^{*}\}$$
where $e=1_H$ is the identity of $S$, the zero element $f$ of
$S$ is the constant zero map and $s=(\alpha_1,u_1)$.  Let $E_e=A$,
$E_f=\{(0,0)\}$, $E_{s^{*}}=E_{s^{*}s}=D_{-1}$ and $E_s=E_{ss^{*}}
=D_1$.  Define
the semigroup action $\beta$ of $S$ as in Proposition
\cite{hivr}.  Then $\beta_s$ is the forward shift, $\beta_{s^{*}}$ is the
backward shift and $\beta_t$ is the identity map for all other
$t\in S$.  As we have seen in Theorem \cite{hi1} the
crossed product $A\times_{\beta}S$ is isomorphic to the matrix
algebra $M_2$.  \eop

\uj Notice that in the last example the semigroup $S$ is
isomorphic to the inverse semigroup generated by
$\{\alpha_n:n\in {\bf N}\}$ as well as to the inverse semigroup generated
by $\{u_n:n\in {\bf N}\}$.  Based upon experience with group actions,
it might seem natural to expect that $S$ is isomorphic to
the inverse semigroup generated by the range of $u$
whenever $\pi\times u$ is a faithful representation of $A\times_{
\alpha}G$.
Perhaps surprisingly this is not the case.  All three
semigroups can be non-isomorphic as the following
example shows.

\example Let $A=C[0,1]$, $G={\bf Z}_2$, $D_0=A$, and let $\alpha_
1$ be the
identity map on $D_1=\{x\in A:x(1)=0\}$.  We construct a
faithful representation $\pi\times u$ of the partial crossed
product $A\times_{\alpha}G$.  Let $\pi$ be the representation of $
A$ on the
Hilbert space $L^2[0,1]\times L^2[0,1]$ defined by $\pi (f)=\left
(\matrix{f&0\cr
0&f\cr}
\right)$ and
let $u_0=\left(\matrix{1&0\cr
0&1\cr}
\right)$ and $u_1=\left(\matrix{0&1\cr
1&0\cr}
\right)$.  By [\cite{mcc},
Propositions 3.4 and 4.2] $\pi\times u$ is faithful since ${\bf Z}_
2$ is
amenable.  The inverse semigroup generated by $\{u_0,u_1\}$ is
isomorphic to ${\bf Z}_2$.  It is clear that $\{\alpha_0,\alpha_1
\}$ is a
semilattice, hence is definitely not isomorphic to the
inverse semigroup $\{u_0,u_1\}$.  The inverse semigroup
generated by $\{(\alpha_0,u_0),(\alpha_1,u_1)\}$ contains three elements
$\{(\alpha_0,u_0),(\alpha_1,u_1),(\alpha_1,u_0)\}$.  \eop

\uj Although every partial crossed product is isomorphic
to a crossed product by an action of a unital inverse
semigroup, this semigroup action may not be unique up
to isomorphism.  For all we know different faithful
representations $\Pi =\pi\times u$ of the crossed product $A\times_{\a}
G$
could generate essentially different semigroup actions.
If we want to talk about a canonical semigroup action
associated with $A\times_{\a}G$ then we can choose $\Pi$ to be the
universal representation of $A\times_{\a}G$.

%\vfil
%\eject
%\double
\uj
\goodbreak
\centerline{\nagy References}
\writeitem{}{References}
%\newsection{References}
\def\i {\advance\sorszam by1 \item{[\the\sorszam]} }
\sorszam=0
\uj
\i R.S. Doran, V.A. Belfi, {\it Characterizations of $C^*$-algebras},
Marcel Dekker, New York, 1986.
\hivatt{coh}
\i R. Exel, {\it Circle actions on $C^*$-algebras, partial automorphisms
 and a generalized Pimsner-Voiculescu exact sequence}, J. Funct. Anal.
{\bf 122 } (1994), 361-401.
\hivatt{ex}
\i J.M. Howie, {\it An introduction to semigroup theory}, Academic Press,
London, 1976.
\hivatt{how}
\i K. McClanahan, \it K-theory for partial crossed products by discrete
 groups, \rm J. Funct. Anal. to appear.
\hivatt{mcc}
\i J.R. Wordingham, {\it The left regular *-representation of an inverse
 semigroup}, Proc. Amer. Math. Soc. {\bf 86} (1982), 55-58.
\hivatt{wordi}
\nobreak
\bigskip
\nobreak
{\it E-mail address: }{\tt nandor.sieben@asu.edu}
\vfil
\eject

\def\redef#1#2{\line{#1 -- #2 \hfil}}
\input cite.inc
\end